\documentclass[12pt]{amsart}

\usepackage{amssymb}

\newtheorem{theorem}{Theorem}[section]
\newtheorem{lemma}[theorem]{Lemma}
\newtheorem{proposition}[theorem]{Proposition}
\newtheorem{definition}[theorem]{Definition}
\newtheorem{remark}[theorem]{Remark}
\newtheorem{corollary}[theorem]{Corollary}
\newtheorem{example}[theorem]{Example}

\numberwithin{equation}{section}

\begin{document}
\baselineskip=16.5pt

\title[Equivariant reduction to Levi subgroup]{On the
equivariant reduction of structure group of
a principal bundle to a Levi subgroup}

\author[I. Biswas]{Indranil Biswas}

\address{School of Mathematics, Tata Institute of Fundamental
Research, Homi Bhabha Road, Mumbai 400005, India}

\email{indranil@math.tifr.res.in}

\email{param@math.tifr.res.in}

\author{A. J. Parameswaran}

\date{}

\begin{abstract}

Let $M$ be an irreducible projective variety
over an algebraically closed field $k$ of
characteristic zero equipped with an action of a group
$\Gamma$. Let $E_G$ be a principal $G$--bundle over $M$,
where $G$ is a connected reductive algebraic group over $k$,
equipped with a lift of the action of $\Gamma$ on $M$. We
give conditions for $E_G$ to admit a $\Gamma$--equivariant
reduction of structure group to $H$, where $H\, \subset\, G$
is a Levi subgroup. We show that for $E_G$, there is a naturally
associated conjugacy class of Levi subgroups of $G$. Given
a Levi subgroup $H$ in this conjugacy class, $E_G$ admits a
$\Gamma$--equivariant reduction of structure group to $H$,
and furthermore, such a reduction is unique up to an automorphism
of $E_G$ that commutes with the action of $\Gamma$.

\end{abstract}

\maketitle

\section{Introduction}

A holomorphic $G$--bundle over ${\mathbb C}{\mathbb P}^1$ 
admits a holomorphic reduction of structure group to a maximal
torus of $G$, where $G$ is a complex reductive group \cite{Gr}.
In particular, any holomorphic vector bundle over
${\mathbb C}{\mathbb P}^1$ splits as a direct sum of line bundles.
If $V$ is a holomorphic vector bundle over ${\mathbb C}{\mathbb P}^1$
equipped with a lift, as vector bundle automorphisms, of the standard
action of the diagonal matrices in
$\text{SL}(2,{\mathbb C})$ on ${\mathbb C}{\mathbb P}^1$,
then $V$ decomposes as a direct sum of holomorphic line subbundles
each left invariant by the action of the torus \cite{Ku}. More
generally, let $E_G$ be a principal $G$--bundle over an irreducible
complex projective variety $M$, where $G$ is a complex reductive group,
with both $E_G$ and $M$ equipped with algebraic actions of a
connected complex algebraic group $\Gamma$; the action of
$\Gamma$ on $E_G$ commutes with the action of $G$ and the projection
of $E_G$ to $M$ is $\Gamma$--equivariant. Assume that
$E_G$ admits a reduction of structure group to a maximal torus $T$
of $G$. Then $E_G$ admits a $\Gamma$--equivariant reduction
of structure group to $T$
if and only if the action of $\Gamma$ on the automorphism group
of $E_G$ leaves a maximal torus invariant \cite{BP}.

The aim here is to investigate conditions under which
a principal $G$--bundle over a projective variety equipped
with an action of a group $\Gamma$ admits a $\Gamma$--equivariant
reduction of structure group to a Levi subgroup of $G$.

Let $M$ be an irreducible projective
variety over an algebraically closed field $k$ of
characteristic zero on which a group $\Gamma$ acts
as algebraic automorphisms. Let $G$ be a connected
reductive linear algebraic group over $k$ and $E_G$
a principal $G$--bundle over $M$. Let $\text{Aut}(E_G)$
denote the group of all automorphisms of $E_G$.
Suppose we are given a lift
of the action of $\Gamma$ on $M$ to $E_G$ that commutes
with the action of $G$. More precisely, the automorphism
of $E_G$ defined by any $\gamma\,\in\, \Gamma$ is an
algebraic automorphism of $G$--bundle over the
automorphism of $M$ defined by $\gamma$. The action
of $\Gamma$ on $E_G$ induces an action on
$\text{Aut}(E_G)$ through group automorphisms.

A torus is a product of copies of the multiplicative
group ${\mathbb G}_m$.
By a Levi subgroup of $G$ we will mean the centralizer of
some torus of $G$.
Let $E_H\, \subset\, E_G$ be a reduction
of structure group to a Levi subgroup $H$. We will denote by
$Z_0(H)$ the connected component of the center of $H$
containing the identity element. So $Z_0(H)$ is contained in
the automorphism group of $E_H$ and hence contained in
$\text{Aut}(E_G)$.

We prove that $E_H$ is left invariant by the action of $\Gamma$
on $E_G$ if and only if $\Gamma$ acts trivially on the subgroup
$Z_0(H)\, \subset\, \text{Aut}(E_G)$ (Theorem \ref{thm.-1}).

A torus of $\text{Aut}(E_G)$ determines a torus, unique up
to inner automorphism, of $G$. The $G$--bundle $E_G$ admits
a $\Gamma$--equivariant reduction of structure group to the
Levi subgroup $H$ if and only if there is torus $T \, \subset
\, \text{Aut}(E_G)$ satisfying the following two conditions:
\begin{enumerate}
\item{} the action of $\Gamma$ on $T$ is trivial;
\item{} there is a subtorus $T'\, \subset\,Z_0(H)\, \subset\, G$
in the conjugacy class of tori of $G$ defined by $T$ such that
the centralizer of $T'$ in $G$ coincides with $H$.
\end{enumerate}
(See Lemma \ref{lem.-1} and Proposition \ref{prop.-1}.)

Let $T_0\, \subset\, G$ be a torus in the conjugacy class
of tori defined by a maximal torus of $\text{Aut}(E_G)^\Gamma$,
where $\text{Aut}(E_G)^\Gamma\, \subset\, \text{Aut}(E_G)$
is the subgroup fixed pointwise by the action of $\Gamma$. The
conjugacy class of $T_0$
does not depend on the choice of the maximal torus. Let $H_0\,
\subset\, G$ be the centralizer of $T_0$.

In Theorem \ref{thm.-2} we prove that $E_G$
admits a $\Gamma$--equivariant
reduction of structure group $E_{H_0}\, \subset\, E_G$
to $H_0$, which is unique in the following sense:
\begin{enumerate}
\item{} for any $\Gamma$--equivariant reduction of
structure group $E'_{H_0}\, \subset\, E_G$ to $H_0$, 
there is an automorphism $\tau\, \in\,\text{Aut}(E_G)^\Gamma$
such that $\tau (E_{H_0}) \, =\, E'_{H_0}$ as subvarieties
of $E_G$; and

\item{} if $E_H$ is a $\Gamma$--equivariant
reduction of structure group to a Levi subgroup $H\, \subset\, G$,
then there is $g\in G$ and $\tau\, \in\,\text{Aut}(E_G)^\Gamma$
such that $g^{-1}H_0g\, \subset\, H$ and
$E_{H_0} g\, \subset\, \tau(E_H)$.
\end{enumerate}

A theorem due to Atiyah says that for an isomorphism of
a vector bundle over $M$ with any direct sum of indecomposable
vector bundles, the direct summands are unique up to a
permutation of the summands. Theorem 4.1 is an equivariant
principal bundle analog of this result of \cite{At}.

\medskip
\noindent
\textbf{Acknowledgments.}\, We thank Vyjayanthi Chari
for a very useful discussion.

\section{Equivariant reduction to a Levi subgroup}

Let $M$ be an irreducible projective
variety over an algebraically closed field $k$ of
characteristic zero. Let $\Gamma$ be a group acting
on the left on $M$. So we have a map
\begin{equation}\label{def.-phi}
\phi\, :\, \Gamma\times M\, \longrightarrow\, M
\end{equation}
such that for any $\gamma\,\in\, \Gamma$, the map defined
by $x\, \longmapsto\, \phi(\gamma\, ,x)$ is an algebraic
automorphism of $M$, and furthermore, $\phi(\gamma_1\,
\phi(\gamma_2\, ,x))\, =\, \phi(\gamma_1\gamma_2\, ,x)$ for
all $\gamma_1, \gamma_2\, \in\, \Gamma$ and $x\,\in\, M$, with
$\phi(e\, ,x)\, =\,x$, where $e\in \Gamma$ is the identity
element.

Let $G$ be a connected reductive linear algebraic group over
the field $k$ and $E_G$ a principal $G$--bundle over $M$. Let
\begin{equation}\label{def.-f}
f\, :\, E_G\, \longrightarrow\, M
\end{equation}
be the projection. We will denote by $\text{Aut}(E_G)$ the
group of all automorphisms of the $G$-bundle $E_G$ (over the
identity automorphism of $M$). So, $\tau (z)g \, =\, \tau (zg)$
and $f(\tau (z)) \, =\, f(z)$ for all $\tau\,\in\,
\text{Aut}(E_G)$ and $z\,\in\, E_G$. Note that
$\text{Aut}(E_G)$ is an affine algebraic group over $k$.
After fixing a faithful representation of $G$, the group
$\text{Aut}(E_G)$ gets identified with a closed subgroup of the
automorphism group of the associated vector bundle; the
automorphism group of a vector bundle is a Zariski open dense
subset
in the vector space defined by the space of all global
endomorphisms of the vector bundle. The group
$\text{Aut}(E_G)$ is,
in fact, the space of all global section of the adjoint bundle
$$
\text{Ad}(E_G)\, :=\, E_G\times^G G \, =\, (E_G\times G)/G
$$
(the action of any $g\in G$ sends any point
$(z\, ,g')\,\in\, E_G\times G$ to $(zg\, , g^{-1}g'g)$). Let
$$
\text{Aut}^0(E_G)\, \subset\, \text{Aut}(E_G)
$$
be the connected component containing the identity element. So
$\text{Aut}^0(E_G)$ is a connected affine algebraic group over $k$.

Assume that $E_G$ is equipped
with a lift of the action of $\Gamma$ on $M$. So the map
\begin{equation}\label{def.-Phi}
\Phi\, :\, \Gamma\times E_G\, \longrightarrow\, E_G
\end{equation}
defining the action has the property that
for any $\gamma\,\in\, \Gamma$ the map defined by
$z\, \longmapsto\, \Phi(\gamma\, ,z)$ is an algebraic
automorphism of $E_G$ that commutes with the action of $G$
on $E_G$, and $f\circ\Phi(\gamma\, ,z) \,=\, \phi(\gamma\, ,f(z))$,
where $f$ is as in \eqref{def.-f}. Note that the action of $\Gamma$
on $E_G$ induces an action of $\Gamma$ on $\text{Aut}(E_G)$
through algebraic group automorphisms. More precisely, the action
of any $\gamma\,\in\, \Gamma$ sends
$F\, \in\, \text{Aut}(E_G)$ to the automorphism defined by
$$
z\, \longmapsto\, \Phi(\gamma\, ,F(\Phi(\gamma^{-1}\, ,z)))\, ,
$$
where $\Phi$ is as in \eqref{def.-Phi}.

The group $\Gamma$ acts on the adjoint bundle
$\text{Ad}(E_G)$ as follows: the action of
any $\gamma\, \in\, \Gamma$ sends $(z\, ,g) \, \in\, E_G\times G$
to
\begin{equation}\label{def.-ac.}
(\Phi(\gamma\, ,z)\, ,g)\,\in\, E_G\times G\, ;
\end{equation}
this descends to an action of $\Gamma$ on the quotient
$\text{Ad}(E_G)\,=\, (E_G\times G)/G$. This descended action
of $\Gamma$ on $\text{Ad}(E_G)$ lifts the action
of $\Gamma$ on $M$, and it clearly preserves the algebraic
group structure of the fibers of $\text{Ad}(E_G)$. The action
of $\Gamma$ on $\text{Ad}(E_G)$ induces an action of
$\Gamma$ on the space of all global sections of $\text{Ad}(E_G)$,
namely $\text{Aut}(E_G)$. It is straight--forward to check that
this induced action on $\text{Aut}(E_G)$ coincides with
the earlier defined action of $\Gamma$ on $\text{Aut}(E_G)$.

Let
\begin{equation}\label{eq.-ga.-fi.}
\text{Aut}(E_G)^\Gamma\, \subset\, \text{Aut}(E_G)
\end{equation}
be the subgroup that is fixed pointwise by the
action of $\Gamma$ on $\text{Aut}(E_G)$. Since the action
of each $\gamma\,\in\, \Gamma$ is an algebraic
automorphism of $\text{Aut}(E_G)$, the subgroup
$\text{Aut}(E_G)^\Gamma$ is Zariski closed.

A reduction of structure group of $E_G$ to a closed
subgroup $H\, \subset\, G$ is a section of $E_G/H$ over $M$,
or equivalently, a closed subvariety $E_H\, \subset\, E_G$
closed under the action of $H$ such that the $H$ action
on $E_H$ defines a principal $H$--bundle over $M$.

\begin{definition}
{\rm A reduction of structure group $E_H\, \subset\, E_G$
to $H$ is called $\Gamma$--\textit{equivariant} if the
subvariety $E_H$ is
left invariant by the action of $\Gamma$ on $E_G$.}
\end{definition}

Since the actions of $\Gamma$ and $G$ on $E_G$ commute, there is an
induced action of $\Gamma$ on $E_G/H$. It is easy to see that $E_H$
is a $\Gamma$--equivariant reduction
of structure group if and only if the section
over $M$ of the bundle $E_G/H\, \longrightarrow\, M$
defined by $E_H$ is fixed by the action of $\Gamma$ on the space
of all sections of $E_G/H$ induced by the action on $E_G/H$.

By a Levi subgroup of $G$ we will mean the
centralizer in $G$ of some torus of $G$. Recall that a torus
is a product of copies of ${\mathbb G}_m$ or the trivial
group. For a Levi subgroup
$H\, \subset\, G$, the centralizer in $G$ of the connected
component of the center of $H$ containing the identity element
coincides with $H$ (see \cite[\S 3]{SS}).
If $H\, \subset\, G$ is a Levi subgroup, then
there is a parabolic subgroup $H\, \subset\, P\, \subset\, G$ such
that $H$ projects isomorphically to the Levi quotient of $P$.
Conversely, if $H$ is a reductive subgroup of a parabolic
subgroup $P\, \subset\, G$ such that $H$ projects isomorphically to
the Levi quotient of $P$, then $H$ is a Levi subgroup of $G$.
Note that if we take the torus to be the trivial group, then
the corresponding Levi subgroup is $G$ itself, and hence in that
case the corresponding parabolic subgroup is $G$.

Take a Levi subgroup $H\, \subset\, G$. Let
$$
Z_0(H)\, \subset\, H
$$
be the connected component of the center of $H$ containing the
identity element. Let
\begin{equation}\label{eq.-red.}
E_H\, \subset\, E_G
\end{equation}
be a reduction of structure group of $E_G$ to $H$. We have
\begin{equation}\label{eq.-incl.}
Z_0(H)\, \subset\, \text{Aut}^0(E_H) \, \subset\,
\text{Aut}^0(E_G)\, ,
\end{equation}
where $\text{Aut}^0(E_H)$ is the connected component of the
group of all automorphisms of the $H$--bundle $E_H$ containing
the identity automorphism; the group $Z_0(H)$ acts on $E_H$
as translations (using the action of $H$ on $E_H$), which
makes $Z_0(H)$ a subgroup of $\text{Aut}^0(E_H)$.

\begin{theorem}\label{thm.-1}
If the reduction $E_H$ in \eqref{eq.-red.} is
$\Gamma$--equivariant, then the subgroup
$Z_0(H)\, \subset\, {\rm Aut}^0(E_G)$ in \eqref{eq.-incl.}
is contained in ${\rm Aut}(E_G)^\Gamma$ (defined in
\eqref{eq.-ga.-fi.}).

Conversely, if $Z_0(H)\, \subset\, {\rm Aut}^0(E_G)\cap
{\rm Aut}(E_G)^\Gamma$, then the reduction
$E_H$ in \eqref{eq.-red.} is $\Gamma$--equivariant.
\end{theorem}

\begin{proof}
Assume that the reduction $E_H$ in \eqref{eq.-red.} is
$\Gamma$--equivariant. For any $\gamma\,\in\, \Gamma$, let
$\Phi_{\gamma}$ be the automorphism of the variety
$E_H$ defined by the action of $\gamma$. The automorphism
$g\, \in\, Z_0(H) \, \subset\, \text{Aut}^0(E_G)$
preserves $E_H$, and on $E_H$ it coincides with the
map $z\, \longmapsto\, zg$. Let $S_g$ be the automorphism of
the $H$--bundle $E_H$
defined by $z\, \longmapsto\, zg$. Since the actions of
$G$ and $\Gamma$ on $E_G$ commute, we have
$$
\Phi_{\gamma}\circ S_g\circ \Phi^{-1}_{\gamma}\, =\,
S_g\circ \Phi_{\gamma}\circ \Phi^{-1}_{\gamma} \, =\, S_g
$$
on $E_H$. Therefore, the two automorphisms, namely
$g\,\in\, \text{Aut}^0(E_G)$
(in \eqref{eq.-incl.}) and the image of $g$ by the action 
$\gamma$ on $\text{Aut}(E_G)$, coincide over $E_H\, \subset\, E_G$.
Consequently, these two automorphisms of $E_G$ coincide. In other
words, the action of $\Gamma$ on ${\rm Aut}(E_G)$
fixes the subgroup $Z_0(H)$ pointwise. This completes the
proof of the first part.

Assume that $\Gamma$ acts trivially
on the subgroup $Z_0(H)\, \subset\,
\text{Aut}(E_G)$ defined in \eqref{eq.-incl.}.
Take a closed point $x\in M$. We will show that the evaluation map
\begin{equation}\label{eq.-ev.}
f_x\, :\, Z_0(H)\, \longrightarrow\, \text{Ad}(E_G)_x
\end{equation}
is injective, where $\text{Ad}(E_G)_x$ is the fiber of
$\text{Ad}(E_G)$ over $x$; the map $f_x$ sends any $s\in Z_0(H)$ to
the evaluation at $x$ of the corresponding section
(as in \eqref{eq.-incl.}) of $\text{Ad}(E_G)$.

To prove that $f_x$ is injective,
fix a finite dimensional faithful $G$--module $V$ over $k$.
Let
$$
E_V\, :=\, (E_G\times V)/G
$$
be the vector bundle over $M$ associated to $E_G$
for the $G$--module $V$; the action of any $g\in G$ sends
$(z\, ,v)\, \in\, E_G\times V$ to $(zg\, ,g^{-1}v)$.
Take any $\sigma\,\in\, Z_0(H)\, \subset\, \text{Aut}^0(E_G)$.
So $\sigma$ gives an automorphism
$$
\sigma'\, \in\, H^0(M,\, \text{Isom}(E_V))
$$
of the vector bundle; the automorphism of $E_G\times V$
that sends any $(z\, ,v)\, \in\, E_G\times V$
to $(\sigma(z)\, ,v)$ descends to an automorphism of $E_V$.

Since $M$ is complete and irreducible, there
are no nonconstant functions on it. Therefore, the
coefficients of the characteristic polynomial of the endomorphism
$$
\sigma'(y) \, \in\, \text{End}((E_V)_y)\, ,
$$
where $y\in M$ is a closed
point, are independent of $y$. Also, since $\sigma$ is
an element of a torus, namely $Z_0(H)$, the
endomorphism $\sigma'(y)$ is semisimple.

If $f_x(\sigma) \, =\, \text{Id}_{(E_G)_x}$, where
$f_x$ is defined in \eqref{eq.-ev.},
then clearly $\sigma'(x)\, =\, \text{Id}_{(E_V)_x}$. Therefore,
in that case, all the eigenvalues of $\sigma'(y)$ are $1$ for
all $y\in M$. Since all $\sigma'(y)$ is semisimple with all
eigenvalues $1$, it follows
immediately that $\sigma'(y)$ is the identity automorphism
of $(E_V)_y$ for each $y\in M$.

Since $V$ is a faithful $G$--module and
$\sigma'$ is the identity automorphism of $E_V$, we conclude that
$\sigma$ is the identity automorphism of $E_G$. This proves
that the homomorphism $f_x$ defined
in \eqref{eq.-ev.} is injective.

Therefore, using the evaluation map, $M\times Z_0(H)\, \subset\,
\text{Ad}(E_G)$ is a subgroup--scheme. Since $Z_0(H)$ is preserved
by the action of $\Gamma$ on $\text{Aut}(E_G)$, it follows
immediately that the action of $\Gamma$ on $\text{Ad}(E_G)$
leaves this subgroup--scheme invariant.

Fix an element
$g_0\, \in\, Z_0(H)$ such that the Zariski closure of the
group generated by $g_0$ coincides with $Z_0(H)$. Since
$H$ is the centralizer of the subgroup $Z_0(H)\, \subset\, G$,
and the algebraic subgroup generated by $g_0$
coincides with $Z_0(H)$, we conclude that $H$ coincides
with the centralizer of $g_0\, \in\, G$.

Let
\begin{equation}\label{ad.-const.}
F\, :\, E_G\times G\, \longrightarrow\, \text{Ad}(E_G)
\, :=\, (E_G\times G)/G
\end{equation}
be the quotient map. Let
\begin{equation}\label{def.-h.-F}
\widehat{F}\, :=\, F^{-1}(\text{image}(\hat{g}_0))\,\subset\,
E_G\times G
\end{equation}
be the subvariety, where
\begin{equation}\label{def.-hat-g}
\hat{g}_0\, :\, M\, \longrightarrow\, \text{Ad}(E_G)
\end{equation}
is the section defined by the above element $g_0$ using the
inclusion $Z_0(H)\, \hookrightarrow\, \text{Aut}(E_G)$ in
\eqref{eq.-incl.}. Set
\begin{equation}\label{wdt.-E}
\widehat{E}\, :=\, \widehat{F}\cap (E_G\times \{g_0\})\, \subset\,
E_G\times G\, ,
\end{equation}
where $\widehat{F}$ is defined in \eqref{def.-h.-F}, and let
\begin{equation}\label{def.-E-p}
E'\, \subset\, E_G
\end{equation}
be the image of $\widehat{E}$ (constructed in
\eqref{wdt.-E}) by the
projection of $E_G\times G$ to $E_G$ defined by
$(z\, ,g)\, \longmapsto\, z$.

Since $\Gamma$ acts trivially on the subgroup
$Z_0(H)\, \hookrightarrow\, \text{Aut}(E_G)$, the image of
the map $\hat{g}_0$ in \eqref{def.-hat-g} is left invariant
by the action of $\Gamma$ on $\text{Ad}(E_G)$. Since the
action of $\Gamma$ on $\text{Ad}(E_G)$ is the
descent, by the projection $F$ in \eqref{ad.-const.}, of the
diagonal action on $E_G\times G$ with
$\Gamma$ acting trivially on $G$, it follows that $E'$ in
\eqref{def.-E-p} is left invariant by the action of
$\Gamma$ on $E_G$.

Since $E'$ is left invariant by $\Gamma$,
the theorem follows
once we show that $E'$ coincides with the subvariety
$E_H$ in \eqref{eq.-red.}.

To prove that $E'\, =\, E_H$, first note that
$$
E_H\times\{g_0\}\, \subset\, \widehat{F}\, \subset\,
E_G\times G
$$
with $\widehat{F}$ defined in \eqref{def.-h.-F}. Indeed,
the automorphism of $E_H$ defined by $g_0$ sends
any $z\, \in\, E_H$ to $zg_0$ (since $g_0$ is in the
center of $H$, this commutes with the
action of $H$ and hence it is an automorphism of $E_H$). This
immediately implies that
$E_H\times \{g_0\}\, \subset\, \widehat{F}$.
Consequently, we have $E_H\, \subset\, E'$. On the other
hand, for any $x\in M$ and $w\,\in\, E'\cap (E_G)_x$
it can be shown that the fiber of $E'$ over $x$
is contained in the orbit of $w$ for the action of
the centralizer of $g_0$ in $G$. Indeed, if
$F(w'\, ,g')\, =\, F(w'g\, ,g')$, where $g, g'\,\in\, G$,
$w'\,\in\, (E_G)_x$ and $F$ as in \eqref{ad.-const.}, then
$gg'g^{-1} \, =\, g'$, this being an immediate consequence
of the definition of $F$. Therefore, if
$w, wg\,\in\, E'$, with $g\in G$, then $g^{-1}g_0g \, =\, g_0$.

We already noted that the centralizer of $g_0$ in $G$
coincides with $H$. We also saw that $E_H\, \subset\, E'$.
Therefore, the above observation that any two points
of $E'$ over a point $x\in M$ differ by an element of
the centralizer of $g_0$ implies that $E_H\, =\, E'$.
This completes the proof of the theorem.
\end{proof}

\begin{example}\label{exa.-1}
{\rm It may happen that $\Gamma$ preserves
the subgroup $Z_0(H)\, \subset\, \text{Aut}^0(E_G)$
in \eqref{eq.-incl.}, but does not preserve $Z_0(H)$
pointwise. We give an example.

Fix a maximal torus $T\, \subset\, G$.
Take $\Gamma$ to be the normalizer $N(T)$ of $T$ in $G$,
and equip $M$ with the trivial action of $\Gamma$;
let $G$ be such that $N(T)\,\not=\, T$. Set $E_G$
to be the trivial $G$--bundle $M\times G$. The group $N(T)$
acts on $M\times G$ as left translations
of $G$. So the induced action of
$N(T)$ on $\text{Aut}(E_G)\, =\, G$ is the conjugation action.
Set $H\, =\, T$. The reduction of structure group of $E_G$
to $T$ defined by the inclusion $M\times T\, \hookrightarrow\,
M\times G$ has the property that the subgroup
$$
Z_0(H) \, =\, T\, \subset\, G\, =\, \text{Aut}(E_G)
$$
(defined in \eqref{eq.-incl.})
is left invariant by the action of $N(T)$ (in this
case it is the adjoint action of $N(T)$
on $G$). However no reduction of structure group
of $E_G$ to $T$ is left invariant by the action $N(T)$.

The automorphism group of a torus is a discrete group. Therefore,
if $\Gamma$ is a connected algebraic group acting algebraically
on $E_G$, then $\Gamma$ acts trivially on $Z_0(H)$ provided
$Z_0(H)$ is preserved by $\Gamma$.}
\end{example}

\begin{proposition}\label{propo0}
Let $T\, \subset\, {\rm Aut}^0(E_G)\cap {\rm Aut}^0(E_G)^\Gamma$
be a torus such that there is an element
$g\,\in\, {\rm Aut}^0(E_G)$ satisfying the condition that
$g^{-1}Tg\, =\, Z_0(H)$, with $Z_0(H)$
constructed in \eqref{eq.-incl.} for the reduction $E_H$
in \eqref{eq.-red.}.
Then $E_G$ admits a $\Gamma$--equivariant
reduction of structure group to the Levi subgroup $H$.
\end{proposition}

\begin{proof}
Take $T$ and $g$ as above. So, the image
$$
E'_H \, :=\, g(E_H)\, \subset\, E_G
$$
is a reduction of structure group of $E_G$ to $H$, where
$E_H$ is the reduction in \eqref{eq.-red.}. Take any automorphism
$\tau$ of the principal
$H$--bundle $E_H$. Using the reduction $E_H$ in
\eqref{eq.-red.}, the automorphism $\tau$ gives
an automorphism $\tau_1$ of the $G$--bundle $E_G$. On the
other hand, using the above reduction $E'_H\, \subset\, E_G$
together with the
isomorphism of $E_H$ with $E'_H$ defined by $z\, \longmapsto\,
g(z)$
the automorphism $\tau$ gives an automorphism $\tau_2$ of $E_G$.
It is easy to see that $\tau_2\, =\, g\tau_1 g^{-1}$.

Therefore, if we substitute $E_H$ by $E'_H$, then the subgroup
$Z_0(H)\, \subset\, \text{Aut}^0(E_G)$ in \eqref{eq.-incl.}
gets replaced by $gZ_0(H)g^{-1}$. Now, the second part of
Theorem \ref{thm.-1} says that $E'_H$ is
left invariant by the action of $\Gamma$ on $E_G$. This
completes the proof of the proposition.
\end{proof}

\section{Levi reduction from tori in $\text{Aut}(E_G)^\Gamma$}

Let $T\,\subset\,\text{Aut}^0(E_G)$ be a torus. From the
proof of Theorem \ref{thm.-1} it
can be deduced that $T$ determines a torus, unique
up to an inner automorphism, in $G$. This will be explained
below with more details.

Fix a point $x\in M$. We saw in the proof of Theorem
\ref{thm.-1} that the evaluation map 
\begin{equation}\label{fx-T}
f_x\, :\, T\, \longrightarrow\, \text{Ad}(E_G)_x
\end{equation}
is injective. Since $\text{Ad}(E_G) \, =\, (E_G\times G)/G$,
if we fix a point $z\in (E_G)_x$, then the quotient map
$F$ (defined in \eqref{ad.-const.}) gives an isomorphism of
$\{z\}\times G$ with $\text{Ad}(E_G)_x$. This identification
of $G$ with $\text{Ad}(E_G)_x$ constructed using $z$ is
an isomorphism of algebraic groups. Furthermore, if we
substitute $z$ by $zg$, $g\in G$, then the corresponding
isomorphism of $G$ with $\text{Ad}(E_G)_x$ is the composition
of the earlier one with the automorphism of
$G$ defined by the conjugation action of $g$. Therefore,
$f_x(T)$, with $f_x$ defined in \eqref{fx-T}, gives a torus
in $G$ up to conjugation.

This torus of $G$, up to conjugation, defined by $f_x(T)$
actually does not depend on the choice of the point $x$. To prove
this, take $z_1, z_2\, \in\, E_G$ with $f(z_i) \, =\, x_i$,
$i=1,2$, where $f$ is as in \eqref{def.-f}. Consider
the evaluation homomorphism
$$
f_{x_i}\, :\, T\, \longrightarrow\, \text{Ad}(E_G)_{x_i}
$$
which is injective. Let
\begin{equation}\label{def.-eq.-h}
h_{z_i}\, :\, T\, \longrightarrow\, G
\end{equation}
be the composition of $f_{x_i}$ with the identification
of $\text{Ad}(E_G)_{x_i}$ with $G$ defined by $z_i$. We want
to show that the two subgroups, namely $\text{image}(h_{z_1})$
and $\text{image}(h_{z_2})$, of $G$ differ
by an inner automorphism of $G$.

Fix a point $t_0\, \in\, T$ such that the Zariski closed
subgroup of $T$ generated by $t_0$ is $T$ itself. For
a finite dimensional $G$--module $V$ over $k$, let
$E_V$ be the vector bundle associated to $E_G$ for $V$ and
$\hat{t}_0$ the automorphism of $E_V$ defined by $t_0
\, \in\, \text{Aut}(E_G)$. From the
definition of the map $h_{z_i}$ it follows
that the automorphism $\hat{t}_0(x_i)$ of $(E_V)_{x_i}$ and the
automorphism of $V$ given by $h_{z_i}(t_0)\, \in\, G$
are intertwined by the isomorphism of $(E_V)_{x_i}$ with
$V$ constructed using $z_i$. (Since $E_V\, =\, (E_G\times V)/G$,
we have an isomorphism of $(E_V)_{x_i}$ with $V$ that sends any
$v\in V$ to the image of $(z_i\, ,v)$.)
We saw in the proof of Theorem \ref{thm.-1} that
the characteristic polynomial of $\hat{t}_0(y)
\, \in\, \text{Isom}((E_V)_{y})$ is independent of $y$.
Therefore, the automorphisms of $V$ defined the two elements
$h_{z_1}(t_0)$ and $h_{z_2}(t_0)$ of $G$ have same characteristic
polynomial.

On the other hand, if $T''\, \subset\, G$ is a maximal
torus, then the algebra of all functions on the affine
variety $T''/W$, where $W\, :=\, N(T'')/T''$ is the Weyl group
with $N(T'')$ the normalizer of $T''$ in $G$,
is generated by trace function of finite dimensional
$G$--modules over $k$ \cite[p. 87, Theorem 2]{St}. Therefore,
$h_{z_1}(t_0)$ and $h_{z_2}(t_0)$ differ by an inner
automorphism of $G$ (since the characteristic polynomials
of $h_{z_1}(t_0)$ and $h_{z_2}(t_0)$ coincide for any
$G$--module). Since $\text{image}(h_{z_i})$ is generated,
as a Zariski closed subgroup, by $h_{z_i}(t_0)$, we conclude
that the two subgroups $\text{image}(h_{z_1})$ and
$\text{image}(h_{z_2})$ differ by an inner automorphism of $G$.

\begin{remark}\label{rem.-3.1}
{\rm Let $E_H\, \subset\, E_G$ be a reduction of structure
group to a Levi subgroup $H\, \subset\, G$. Consider the torus
$Z_0(H)\, \subset\, \text{Aut}^0(E_G)$ in \eqref{eq.-incl.}
corresponding to the reduction $E_H$. By substituting
a point of $E_H$ for the point $z_i$ in \eqref{def.-eq.-h}
we conclude that the map in \eqref{def.-eq.-h} sends
any $g\,\in\, Z_0(H) \, \subset\, H$ to the point
$g\, \in\, \text{Aut}^0(E_G)$ (in terms of
\eqref{eq.-incl.}). Consequently, the torus $Z_0(H)\, \subset\,
G$ is in the conjugacy class of tori given by the torus
$Z_0(H)\, \subset\, \text{Aut}^0(E_G)$ in \eqref{eq.-incl.}.}
\end{remark}

We have the following lemma:

\begin{lemma}\label{lem.-1}
If the $G$--bundle $E_G$ admits a $\Gamma$--equivariant
reduction of structure group to a Levi
subgroup $H\, \subset\, G$, then there is a torus $T\, \subset\,
{\rm Aut}^0(E_G)\cap {\rm Aut}(E_G)^\Gamma$ that
satisfies the condition that $Z_0(H)$
is the torus in $G$ defined, up to conjugation, by $T$.
\end{lemma}

\begin{proof}
Let $E_H\, \subset\, E_G$ be a $\Gamma$--equivariant reduction of
structure group to $H$. The image of $Z_0(H)$ in $\text{Aut}^0(G)$
by the inclusion map in \eqref{eq.-incl.} will be denoted by $T$.
The first part of Theorem \ref{thm.-1} says that $T
\, \subset\,\text{Aut}(E_G)^\Gamma$.

Fix a point $z\, \in\, E_H\,\subset\, E_G$. It is easy to see
that the torus $h_{z}(T)\, \subset\, G$ coincides with $Z_0(H)$,
where $h_{z}$ is defined as in \eqref{def.-eq.-h} (by
composing the evaluation map $T\, \longrightarrow\,
\text{Ad}(E_G)_{f(z)}$, where $f$ is defined in \eqref{def.-f},
with the isomorphism $\text{Ad}(E_G)_{f(z)} \, \longrightarrow\,
G$ defined by $z$). This completes the proof of the lemma.
\end{proof}

In the converse direction we have:

\begin{proposition}\label{prop.-1}
Let $T'\, \subset\, G$ be a torus in the conjugacy class of tori
determined by a torus $T\, \subset\, {\rm Aut}^0(E_G)^\Gamma$ and
$H$ the centralizer of $T'$ in $G$. Then $E_G$ admits a 
$\Gamma$--equivariant reduction of structure group to the Levi
subgroup $H$.
\end{proposition}

\begin{proof}
Fix any point $z\, \in\, E_G$ and consider the homomorphism
\begin{equation}\label{def.-hz}
h_{z}\, :\, T\, \longrightarrow\, G
\end{equation}
as in \eqref{def.-eq.-h}, namely it is the composition of the
evaluation map with the identification, constructed
using $z$, of $G$ with $\text{Ad}(E_G)_{f(z)}$,
where $f$ is defined in \eqref{def.-f}. There is an element
$g\in G$ with $gh_z(T)g^{-1} \, =\, T'$, where $T'$
is as in the statement of the proposition.

Let
\begin{equation}\label{def.-H-pr.}
H_z\, \subset\, G
\end{equation}
be the centralizer of $h_z(T)$,
with $h_z$ defined in \eqref{def.-hz}.
Since $gh_z(T)g^{-1} \, =\, T'$, and the centralizer
of $T'\,\subset\, G$ is $H$, we conclude that
\begin{equation}\label{H-pr.-id.}
gH_zg^{-1} \, =\, H\, .
\end{equation}

Fix an element $t_0\, \in\, T$ such that the Zariski
closure in $T$ of the subgroup generated by $t_0$ is
$T$ itself. As in \eqref{def.-hat-g}, let
$$
\hat{t}_0\, :\, M\, \longrightarrow\, \text{Ad}(E_G)
$$
be the section defined by the automorphism $t_0\,\in\,
\text{Aut}(E_G)$. As in \eqref{def.-h.-F}, set
$$
\widehat{F}\, :=\, F^{-1}(\text{image}(\hat{t}_0))\,\subset\,
E_G\times G\, ,
$$
where $F$ is the projection in \eqref{ad.-const.}. As
in \eqref{wdt.-E}, define
$$
\widehat{E}\, :=\, \widehat{F}\cap (E_G\times \{h_z(t_0)\})
\, \subset\, E_G\times G\, ,
$$
where $h_z$ is defined in \eqref{def.-hz}. Let
\begin{equation}\label{def.-2p}
E'\, \subset\, E_G
\end{equation}
be the image of $\widehat{E}$ by the projection of $E_G\times G$
to $E_G$ defined by $(y\, ,\nu)\, \longmapsto\, y$.

We will show that $E'$ constructed in \eqref{def.-2p} is a
$\Gamma$--equivariant reduction of structure group of $E_G$
to the subgroup $H_z$ defined in \eqref{def.-H-pr.}.

For this, we will first show that $E'$ is closed under the action
of $H_z$ (for the action of $G$ on $E_G$). Note that $h_z(t_0)$
is in the center of $H_z$ (as $H_z$ is the centralizer of $h_z(T)$).
Therefore, for the projection $F$ in \eqref{ad.-const.} we have
\begin{equation}\label{eq.-id.-p.}
F(z_1\, ,h_z(t_0))\, =\, F(z_1g_1\, ,h_z(t_0))
\end{equation}
for all $z_1\in E_G$ and and $g_1\in H_z$. Indeed, the map
$F$ clearly has the property that for
$g, g'\,\in\, G$ and $w'\,\in\, E_G$
\begin{equation}\label{eq.-p2.}
F(w'\, ,g')\, =\, F(w'g\, ,g')
\end{equation}
if and only if $gg'g^{-1} \, =\, g'$. From \eqref{eq.-id.-p.} it
follows immediately
that $E'$ is closed under the action of $H_z$.

It also follows from \eqref{eq.-p2.} that for any point
$y\in M$, the centralizer of $h_z(t_0)$ (in $G$) acts
transitively on the fiber of $E'$ over $y$. Note that since the
Zariski closure of the group generated by $t_0$ is $T$, and
$H_z$ is the centralizer (in $G$) of $h_z(T)$, it follows
immediately that the centralizer of $h_z(t_0)$ is $H_z$.

We still need to show that the fiber of $E'$ over each 
point $y\,\in\, M$ is nonempty. For this note that there
is a point $z'\, \in\, f^{-1}(y)$, with $f$ defined in
\eqref{def.-f}, such that the corresponding homomorphism
$$
h_{z'}\, :\, T\, \longrightarrow\, G
$$
defined as in \eqref{def.-hz} by replacing $z$ by $z'$
has the property that $h_{z'}(t_0) \, =\, h_{z}(t_0)$.
Indeed, this follows from the combination of the
fact that the conjugacy
class of the torus $h_z(T)\, \subset\, G$ is independent
of the choice of the point $z\in E_G$ and the observation
that the two homomorphisms $h_{z_1}$ and $h_{z_1g_1}$ from
$T$ to $G$, where $z_1\,\in\, E_G$ and
$g_1\in G$, differ by the inner automorphism of $G$ defined
by $g_1$. The identity $h_{z'}(t_0) \, =\, h_{z}(t_0)$
immediately implies that $z'$ is in the fiber of $E'$ over $y$.

Consequently, $E'\, \subset\, E_G$ constructed in \eqref{def.-2p}
is a reduction of structure group to $H_z$.

Since the action of
$\Gamma$ on $\text{Aut}(E_G)$ fixes $t_0$, the action of $\Gamma$
on $E_G$ leaves $E'$ invariant.

Finally, from \eqref{H-pr.-id.} it follows immediately
that $E'g^{-1}\, \subset\, E_G$ is a reduction of
structure group to $H$. As $E'$ is left invariant by the
action of $\Gamma$ on $E_G$, and the actions of $\Gamma$ and $G$
on $E_G$ commute, the subvariety $E'g^{-1}\, \subset\,
E_G$ is also left invariant by the action of $\Gamma$. This
completes the proof of the proposition.
\end{proof}

\section{A canonical equivariant Levi reduction}

Let $T\, \subset\, \text{Aut}(E_G)^\Gamma$ be a
connected maximal torus, where $\text{Aut}(E_G)^\Gamma$
is defined in \eqref{eq.-ga.-fi.}. So $T$ is a torus
of $\text{Aut}^0(E_G)$.
We saw in the previous section that $T$ determines a
torus, unique up to an inner conjugation, in $G$. We will
show that this torus
in $G$ (up to conjugation) does not depend on the
choice of the maximal torus $T$.

To prove this, first note that any two
maximal tori of $\text{Aut}(E_G)^\Gamma$ differ by an
inner automorphism of $\text{Aut}(E_G)^\Gamma$.
Consider the maximal torus $g_0Tg_0^{-1}$,
where $g_0\,\in\, \text{Aut}(E_G)^\Gamma$, and fix a point
$z\in E_G$. The point $z$ defines two injective homomorphisms
$$
h_z\, :\, T\, \longrightarrow\, G
$$
and
$$
h'_z\, :\, g_0Tg_0^{-1} \, \longrightarrow\, G
$$
defined as in \eqref{def.-eq.-h} using the
evaluation map and the isomorphism of groups
$$
\phi_z\, :\, \text{Ad}(E_G)_{f(z)}\, \longrightarrow\, G
$$
constructed using $z$, where $f$ is defined in \eqref{def.-f}. From
the construction of $h_z$ and $h'_z$ it follows immediately that
$$
\phi_z(g_0(f(z))) h_z(T)(\phi_z(g_0(f(z))))^{-1} \, =\, 
h'_z(g_0Tg_0^{-1})\, .
$$
Therefore, $h_z(T)$ and $h'_z(g_0Tg_0^{-1})$ differ by an
inner automorphism of $G$. Consequently, the torus of
$G$ determined by a maximal torus of $\text{Aut}(E_G)^\Gamma$
does not depend on the choice of the maximal torus.

Fix a torus $T_0\, \subset\, G$ in the conjugacy class of
tori given by a maximal torus in $\text{Aut}(E_G)^\Gamma$. The
centralizer of $T_0$ in $G$ is a Levi subgroup. This Levi
subgroup of $G$ will be denoted by $H_0$.

\begin{theorem}\label{thm.-2}
The $G$--bundle $E_G$ admits a $\Gamma$--equivariant
reduction of structure group
to the Levi subgroup $H_0$ defined above.

If $H\,\subsetneq \, H_0$ is a Levi subgroup of $G$
properly contained in $H_0$, then $E_G$ does not admit
any $\Gamma$--equivariant reduction of structure group
to $H$.

If $H\, \subset\, G$ is a Levi subgroup such that
$E_G$ admits a $\Gamma$--equivariant reduction of structure group
to $H$, but $E_G$ does not admit
a $\Gamma$--equivariant reduction of structure group to any
Levi subgroup properly contained in $H$, then $H$ is conjugate
to the above defined subgroup $H_0\, \subset\, G$.

If $E_{H_0}\, \subset\, E_G$ and $E'_{H_0}\, \subset\, E_G$
are two $\Gamma$--equivariant reductions of structure group
to $H_0$, then there is an automorphism $\tau\, \in\,
{\rm Aut}(E_G)^\Gamma$ of $E_G$ such that $\tau(E_{H_0})\, =\,
E'_{H_0}\, \subset\, E_G$.
\end{theorem}

\begin{proof}
That $E_G$ admits a $\Gamma$--equivariant
reduction of structure group to $H_0$
follows from the construction in Proposition \ref{prop.-1}.
Fix a maximal torus $T\, \subset\, \text{Aut}^0(E_G)^\Gamma$ and
a point $t_0\,\in\, T$ such that the Zariski closure of the
subgroup of $T$ generated by $t_0$
coincides with $T$. Let $t'_0\,\in\, T_0$ be the
element corresponding to $t_0$ by an isomorphism of $T$
with $T_0$ constructed using an element of $E_G$.
As in \eqref{wdt.-E}, consider
$$
\widehat{E}\, :=\, F^{-1}(\text{image}(\hat{t}_0))\cap
(E_G\times \{t'_0\}) \, \subset\, E_G\times G\, ,
$$
where $F$ is defined in \eqref{ad.-const.} and $\hat{t}_0$
is the section of $\text{Ad}(E_G)$ defined by $t_0$. Finally
the image of $\widehat{E}$ by the projection of $E_G\times G$
to $E_G$ gives a reduction of structure group of $E_G$
to $H_0$. See the proof of Proposition \ref{prop.-1} for the
details.

To prove the second statement, let $H\, \subsetneq H_0$ 
be a Levi subgroup of $G$ properly contained in $H_0$.
So $\dim Z_0(H) \, >\, \dim T_0$, where $Z_0(H)$ is the
connected component of the center of $H$ containing the
identity element (note that $T_0$ is contained in the center
of the bigger Levi subgroup). The first statement in
Theorem \ref{thm.-1} says that if $E_H\, \subset\, E_G$ is
a $\Gamma$--equivariant reduction of structure group to $H$,
then $\text{Aut}(E_G)^\Gamma$ contains a torus isomorphic
to $Z_0(H)$. This is impossible, since a smaller
dimensional torus, namely $T_0$, is isomorphic to the
maximal torus $T$ and any two maximal tori are isomorphic.

Let $H\, \subset\, G$ be a Levi subgroup as in the third
statement, and let $E_H\, \subset\, E_G$ be a
$\Gamma$--equivariant reduction of structure group to $H$.
The condition on $H$ implies that the torus
$Z_0(H)\, \subset\,\text{Aut}(E_G)^\Gamma$ in \eqref{eq.-incl.}
for the reduction $E_H$ is a maximal torus
of $\text{Aut}(E_G)^\Gamma$. Indeed, that $Z_0(H)\,
\subset\,\text{Aut}(E_G)^\Gamma$ follows from Theorem
\ref{thm.-1}. That $Z_0(H)$ is a maximal torus of
$\text{Aut}(E_G)^\Gamma$ can be seen as follows.
If $T''\, \subset\, \text{Aut}(E_G)^\Gamma$
is a torus with $Z_0(H)\, \subsetneq \, T''$, then
take a torus $T''_1$ in the conjugacy
class of tori of $G$ given by $T''$ such that
$Z_0(H)\, \subset \, T''_1\, \subset\, G$. Let
$H''\, \subset\, G$ be the centralizer of
$T''_1$. Since $Z_0(H)$ is the connected component
of the center of $H$ containing the identity element
and $Z_0(H)\, \subsetneq \, T''_1$ is a proper
subtorus, we conclude that $H'' \, \subsetneq \, H$.
Proposition \ref{prop.-1} says that $E_G$
admits a $\Gamma$--equivariant reduction of structure group
to $H''$. Since $H''$ is a Levi subgroup properly contained
in $H$, this
contradicts the given condition on $H$. Therefore,
$Z_0(H)\, \subset\,\text{Aut}(E_G)^\Gamma$ is a maximal torus.

Since $T_0$, by definition, is in the
conjugacy class of tori of $G$ given by a maximal
torus of $\text{Aut}(E_G)^\Gamma$, using Remark \ref{rem.-3.1}
we conclude that the two tori $T_0$ and $Z_0(H)$ of $G$
are conjugate. Consequently, $H$ and $H_0$ differ by an inner
automorphism of $G$.

Let $E_{H_0}$ and $E'_{H_0}$ be as in the fourth statement.
Consider the inclusion in \eqref{eq.-incl.}. Let $T_1$
(respectively, $T'_1$) be the
image of $T_0$ in $\text{Aut}(E_G)^\Gamma$ for the reduction
$E_{H_0}$ (respectively, $E'_{H_0}$) by \eqref{eq.-incl.}. From
dimension consideration we know that both $T_1$ and
$T'_1$ are maximal tori in $\text{Aut}(E_G)^\Gamma$. Take an
element $\tau\,\in\, \text{Aut}(E_G)^\Gamma$ such that
\begin{equation}\label{eq.-conj.}
T'_1\, =\, \tau^{-1}T_1\tau\, .
\end{equation}

Let $E_H\, \subset\, E_G$ be a $\Gamma$--equivariant
reduction of structure group to a Levi subgroup $H\, \subset\,
G$ and $g_0\, \in\, Z_0(H)$ an element in the connected
component of the center of $H$ containing the identity element
such that $g_0$ generates $Z_0(H)$ as a Zariski closed subgroup
of $G$. In the proof of Theorem \ref{thm.-1} we gave a
reconstruction of $E_H$ from $g_0$ and its image
in $\text{Aut}^0(E_G)$ by \eqref{eq.-incl.}. (In the notation
of the proof of Theorem \ref{thm.-1}, $E'\, \subset\,
E_G$ was constructed in
\eqref{def.-E-p} using $g_0$ and its image in
$\text{Aut}^0(E_G)$, and it was shown there that $E_H$ coincides
with $E'$.)

Fix an element $g_0\, \in\, T_0$ such that Zariski closed subgroup
generated by $g_0$ coincides with $T_0$. Let $g_1$
be the image of $g_0$ in $T_1$ for the above isomorphism of
$T_0$ with $T_1$ constructed using $E_H$. Set
$$
g'_1\, =\, \tau^{-1}g_1\tau\, \in\, T'_1\, ,
$$
where $\tau$ is as in \eqref{eq.-conj.}. Let $g'_0$
be the image of $g'_1$ in $T_0$ for the above isomorphism of
$T'_1$ with $T_0$ constructed using $E'_H$. Following the
construction of $E'$ in \eqref{def.-E-p}, we can
reconstruct $E_H$ (respectively, $E'_H$)
using the pair $(g_0\, , g_1)$ (respectively,
$(g'_0\, , g'_1)$). Using this reconstruction it is
easy to see that
$$E'_H \, =\, \tau^{-1}(E_H)\, ,
$$
where $\tau$ is as in \eqref{eq.-conj.}. This
completes the proof of the theorem.
\end{proof}

\begin{remark}
{\rm If we set $G\, =\, \text{GL}(n,k)$ and $\Gamma\, =\,
\{e\}$, then Theorem \ref{thm.-2} becomes the following
theorem proved in \cite{At}: any vector bundle $V$ over $M$
is isomorphic to a direct sum $\bigoplus_{i=1}^k U_i$
of indecomposable vector bundles, and if
$$
V\, \cong\, \bigoplus_{j=1}^l W_j\, ,
$$
where each $W_j$ is
indecomposable, then $k\,=\,l$ and the collection of
vector bundles $\{W_j\}$ is a permutation of $\{U_i\}$.}
\end{remark}

\begin{remark}
{\rm Let $E_*$ be a parabolic  $G$--bundle over an irreducible
smooth projective variety $X$. Corresponding to $E_*$, there is
an irreducible smooth projective variety $Y$, a finite subgroup
$\Gamma\, \subset\, \text{Aut}(Y)$ with $X\, =\,Y/\Gamma$, and
a principal $G$--bundle $E_G$ over $Y$ equipped with a lift
of the action of $\Gamma$. More precisely, there is a bijective
correspondence between parabolic $G$--bundles and $G$--bundles
with a finite group action on a (ramified) covering (see
\cite{BBN} for the details). Therefore, Theorem \ref{thm.-2}
gives a natural reduction of structure group of a parabolic
$G$--bundle to a Levi subgroup of $G$. This Levi reduction
satisfies all the analogous conditions in Theorem \ref{thm.-2}.}
\end{remark}

\section{The Levi quotient of the automorphism group}

In this final section, we will assume $\Gamma$ to
be a connected algebraic group. We will also assume the action of
$\Gamma$ on $E_G$ to be algebraic, that is,
the map $\phi$ in \eqref{def.-phi} is algebraic;
consequently, the action of $\Gamma$ on $M$ is also algebraic.
Since $\Gamma$ is connected, the action of $\Gamma$ on
$\text{Aut}(E_G)$ preserves the subgroup
$\text{Aut}^0(E_G)$.

Let $U\text{Aut}^0(E_G)$ be the unipotent radical of
the algebraic group $\text{Aut}^0(E_G)$
\cite[p. 125]{Hu}. So the Levi quotient
\begin{equation}\label{def.-Le-a}
L\text{Aut}^0(E_G)\, :=\,
\text{Aut}^0(E_G)/U\text{Aut}^0(E_G)
\end{equation}
is a connected reductive algebraic group over $k$. Let
\begin{equation}\label{def.-psi}
\psi\, :\, \text{Aut}^0(E_G)\, \longrightarrow\,
L\text{Aut}^0(E_G)
\end{equation}
be the quotient map.

{}From the uniqueness of a unipotent radical
it follows immediately
that the action of $\Gamma$ on $\text{Aut}^0(E_G)$ preserves
the subgroup $U\text{Aut}^0(E_G)$. Therefore, we have an induced
action of $\Gamma$ on $L\text{Aut}^0(E_G)$.

Let $\widehat{T}_0\, \subset\, G$ be a torus in the conjugacy
class of tori of $G$ given by a maximal torus in
$\text{Aut}^0(E_G)$. Since any two maximal tori are conjugate,
the conjugacy class of $\widehat{T}_0$ does not depend on the
choice of the maximal torus. Let $\widehat{H}_0$ be the
centralizer of $T_0$ in $G$. Setting $\Gamma\,=\,\{e\}$ in
Proposition we conclude that $E_G$ admits a reduction
of structure group to $\widehat{H}_0$.

\begin{proposition}\label{prop.-2}
If $E_G$ admits a $\Gamma$--equivariant reduction
of structure group to the
Levi subgroup $\widehat{H}_0$, then the induced action of $\Gamma$ on
$L{\rm Aut}^0(E_G)$ (defined in \eqref{def.-Le-a}) factors
through an action of a torus quotient of $\Gamma$.

If $\Gamma$ is reductive and the induced action of $\Gamma$ on
$L{\rm Aut}^0(E_G)$ factors through an action
of a torus quotient of $\Gamma$, then $E_G$ admits a
$\Gamma$--equivariant reduction of structure group
to the Levi subgroup $\widehat{H}_0$.
\end{proposition}

\begin{proof}
Assume that $E_G$ admits a $\Gamma$--equivariant reduction
of structure group to $\widehat{H}_0$. Theorem \ref{thm.-2} says
that there is a maximal torus
$$
T_0\, \subset\, \text{Aut}^0(E_G)
$$
which is left invariant by the action of $\Gamma$
on $\text{Aut}(E_G)$. Consider $\psi(T_0)$, with
$\psi$ defined in \eqref{def.-psi},
which is a maximal torus in
$L\text{Aut}^0(E_G)$. Note that $\psi(T_0)$
is left invariant by the induced action of
$\Gamma$ on $L\text{Aut}^0(E_G)$, as $T_0$
is $\Gamma$--invariant.

Let $ZL\text{Aut}^0(E_G)\, \subset\, L\text{Aut}^0(E_G)$
be the center and
$$
PL\text{Aut}^0(E_G)\, :=\,
L\text{Aut}^0(E_G)/ZL\text{Aut}^0(E_G)
$$
the corresponding adjoint
group. All the automorphisms of $L\text{Aut}^0(E_G)$
connected to the identity automorphism are parametrized
by $PL\text{Aut}^0(E_G)$, with $PL\text{Aut}^0(E_G)$ acting
on $L\text{Aut}^0(E_G)$ as conjugations.

Since $\Gamma$ is connected, we have a homomorphism of
algebraic groups
$$
\rho\, :\, \Gamma\, \longrightarrow\, PL\text{Aut}^0(E_G)
$$
such that the
action of any $g\,\in\, \Gamma$ on $L\text{Aut}^0(E_G)$
is conjugation by $\rho(g)$. Since the action of $\Gamma$ preserves
the maximal torus $\psi(T_0)\, \subset\, L\text{Aut}^0(E_G)$,
and $q\circ\psi(T_0)$ is a maximal torus in
$PL\text{Aut}^0(E_G)$, where
$$
q\, :\, L\text{Aut}^0(E_G)\, \longrightarrow\,
PL\text{Aut}^0(E_G)
$$
is the projection, we conclude that 
$\rho(\Gamma)\, \subset\, q\circ\psi(T_0)$,
where $\psi$ is defined in \eqref{def.-psi}. (The maximal
torus $q\circ\psi(T_0)$ is a finite index subgroup of its
normalizer in $PL\text{Aut}^0(E_G)$.) Therefore, the action
of $\Gamma$ on $L\text{Aut}^0(E_G)$ factors through
the conjugation action  of the torus $\rho(\Gamma)$.

To prove the second statement in the proposition, assume that
the induced action of $\Gamma$ on $L\text{Aut}^0(E_G)$ factors
through the torus quotient $\Gamma\, \longrightarrow\, T_\Gamma$.
We will first
show that the action of $T_\Gamma$ on $L\text{Aut}^0(E_G)$
preserves a maximal torus.

Construct the semi-direct product
$L\text{Aut}^0(E_G)\rtimes T_\Gamma$ using the induced
action of $T_\Gamma$ on $L\text{Aut}^0(E_G)$.
We take a maximal torus
$$
\widehat{T}\, \subset \, L\text{Aut}^0(E_G)\rtimes T_\Gamma
$$
containing $T_\Gamma$ (note that
$T_\Gamma$ is naturally a subgroup of
$L\text{Aut}^0(E_G)\rtimes T_\Gamma$). Finally,
consider the intersection
$$
T_1\, :=\, \widehat{T}\cap L\text{Aut}^0(E_G)
$$
(note that $L\text{Aut}^0(E_G)$ is a normal subgroup of
$L\text{Aut}^0(E_G)\rtimes T_\Gamma$). From its
construction it is immediate that $T_1$ is a maximal
torus of $L\text{Aut}^0(E_G)$ and $T_1$ is left invariant
by the action of $T_\Gamma$ on $L\text{Aut}^0(E_G)$.

Consider
$$
G'\, :=\, \psi^{-1}(T_1)\, \subset\, \text{Aut}^0(E_G)\, ,
$$
where $\psi$ is the projection in \eqref{def.-psi}. Since
$\Gamma$ preserves $T_1\, \subset\, L\text{Aut}^0(E_G)$,
the action of $\Gamma$ on $\text{Aut}^0(E_G)$ preserves
the subgroup $G'$ defined above. Note that $G'$
fits in an exact sequence
$$
e\, \longrightarrow\, U\text{Aut}^0(E_G) \, \longrightarrow\,
G' \, \longrightarrow\, T_1 \, \longrightarrow\, e\, ,
$$
where $U\text{Aut}^0(E_G)$, as before, is the unipotent radical.

A maximal torus of $G'$ is a maximal torus of
$\text{Aut}^0(E_G)$, and since $\Gamma$ is
connected, an algebraic action of $\Gamma$ on a
torus through automorphisms is trivial. Therefore, in view of
Proposition \ref{prop.-1},
to prove the second statement in the proposition it suffices
to show that $\Gamma$ preserves some maximal torus in $G'$.

Denote by ${\mathfrak g}'$ the Lie algebra of $G'$. The action of
$\Gamma$ on $G'$ induces an action of $\Gamma$ on ${\mathfrak g}'$.
Let $\mathfrak u$ (respectively, ${\mathfrak t}_1$)
be the Lie algebra of $U\text{Aut}^0(E_G)$ (respectively, $T_1$).
So the above exact sequence of groups give an exact sequence
\begin{equation}\label{ex.-seq.}
0 \, \longrightarrow\, {\mathfrak u} \, \longrightarrow\,
{\mathfrak g}' \, \stackrel{\beta}{\longrightarrow}\,
{\mathfrak t}_1 \, \longrightarrow\, 0
\end{equation}
of Lie algebras.

Let
$$
{\mathcal V}\, \subset\, {\mathfrak g}'
$$
be the subspace on which $\Gamma$ acts trivially. Note that
${\mathcal V}$ is a Lie subalgebra. The action of
$\Gamma$ on $T_1$ is trivial (as the automorphism group of
$T_1$ is discrete and $\Gamma$ is connected). Therefore,
the induced action of $\Gamma$ on ${\mathfrak t}_1$ is trivial.

Since $\Gamma$ is reductive, any exact sequence of
finite dimensional $\Gamma$--modules over $k$ splits,
in particular, \eqref{ex.-seq.} splits. Since
${\mathfrak t}_1$ is the trivial $\Gamma$--module, we conclude
that the restriction to
the subalgebra ${\mathcal V}\, \subset\, {\mathfrak g}'$
of the projection $\beta$ in \eqref{ex.-seq.} is surjective.

Let $G_2\, \subset\, G'$ be the Zariski closed subgroup
generated by the subalgebra ${\mathcal V}$. Since $\Gamma$
acts trivially on $\mathcal V$ we conclude that $G_2$ is
fixed pointwise by the action of $\Gamma$ on $G'$.

Since the projection of ${\mathcal V}$ to ${\mathfrak t}_1$
(by $\beta$ in \eqref{ex.-seq.}) is
surjective, the subgroup $G_2$ projects surjectively to $T_1$.
Take any maximal torus $T_2\, \subset\, G_2$. Since
the projection of $G_2$ to $T_1$ is surjective and the kernel
of the projection $G' \, \longrightarrow\, T_1$ is a unipotent
group, we conclude that $T_2$ is a maximal torus of $G'$.

In other words, $T_2$ is a $\Gamma$--invariant maximal torus
of $G'$. Since a maximal torus in $G'$ is a maximal torus in
$\text{Aut}^0(E_G)$, Proposition \ref{prop.-1} completes the
proof of the proposition.
\end{proof}

It is easy to construct examples showing that the second statement
in Proposition \ref{prop.-2} is not valid for arbitrary
connected algebraic group $\Gamma$.

Proposition \ref{prop.-2} has the following corollary:

\begin{corollary}
If $\Gamma$ does not have a nontrivial torus quotient
(for example, if it is unipotent or semisimple),
and the action of $\Gamma$ on $L{\rm Aut}^0(E_G)$ is
nontrivial, then $E_G$ does not admit any
$\Gamma$--equivariant reduction of structure group
to $\widehat{H}_0$, provided $\widehat{H}_0\, \not=\, G$.
\end{corollary}


\end{document}